\documentclass[12pt]{amsart}
\usepackage{amssymb,upref,enumerate}

\usepackage{amsmath,amssymb,amscd,amsthm,amsxtra}
\usepackage{latexsym}
\usepackage[dvips]{graphics,epsfig}
\usepackage{color}

\newcommand{\com}{\mathbb{C}}
\newcommand{\na}{\mathbb{N}}
\newcommand{\re}{\mathbb{R}}
\newcommand{\ent}{\mathbb{Z}}

\newcommand{\fr}[2]{{\textstyle \frac{#1}{#2}}}

\newcommand{\norm}[2]{\left\|#1\right\|_{#2}}

\newcommand{\rn}{{\mathbb R}^n}

\newcommand{\zn}{\mathbb Z^n}
\newcommand{\zz}{\mathbb Z}

\newcommand{\sw}{{\mathcal{S}}(\rn)}

\newcommand{\eps}{\varepsilon}

\newcommand{\fh}{\widehat{f}}
\newcommand{\gh}{\widehat{g}}

\newcommand{\Bzinf}{\dot{B}^{0,\infty}_{\infty}}

\newtheorem{thm}{Theorem}

\theoremstyle{remark}

\newtheorem{remark}{Remark}
\numberwithin{equation}{section}

\begin{document}

\title[On an endpoint Kato-Ponce inequality] {A remark on an endpoint Kato-Ponce inequality}
\author{Loukas Grafakos, Diego Maldonado, \and Virginia Naibo}

\address{Department of Mathematics, University of Missouri, Columbia, MO 65211, USA}
\email{grafakosl@missouri.edu}

\address{Department of Mathematics, Kansas State University, Manhattan, KS 66506, USA}
\email{dmaldona@math.ksu.edu}

\address{Department of Mathematics, Kansas State University, Manhattan, KS 66506, USA}
\email{vnaibo@math.ksu.edu}

\thanks{Third author partially supported by NSF under grant DMS 1101327}


\date{\today}


\begin{abstract}
This note introduces bilinear estimates intended as a step towards an $L^\infty$-endpoint Kato-Ponce inequality. In particular, a bilinear version of the classical Gagliardo-Nirenberg interpolation inequalities for a product of functions is proved. 
\end{abstract}

\maketitle

\section{Introduction and main result}

The following inequality appears to be missing from the vast literature on a class of inequalities known as Kato-Ponce inequalities or fractional Leibniz rules: For every  $s>0$ there exists  $C>0$, depending only on $s$ and dimension $n$, such that 
\begin{equation}\label{KPinfty}
\norm{D^s(fg)}{L^\infty}  \leq C \left( \norm{D^sf}{L^\infty} \norm{g}{L^\infty} + \norm{D^sg}{L^\infty} \norm{f}{L^\infty}\right), \quad \text{for all } f, g \in \sw,
\end{equation}
where $D^s$ is the $s$-derivative  operator\footnote{The notation $D^s$ seems to be standard for this operator although other notations include $|D|^s,$   $|\nabla|^s$ and $(-\Delta)^{\frac{s}{2}}.$} defined for $h \in \sw$ as
$$
\widehat{D^sh}(\xi):= |\xi|^s \hat{h}(\xi), \quad  \forall \xi \in \rn.
$$ 
Inequality \eqref{KPinfty} represents an endpoint case of inequalities of Kato-Ponce type (see \cite{BBR, BMMN, CWein, CRT, GO, GK, KaPo, KPVe, KoTa1} and references therein) and we do not know whether it holds true or not. Moreover, the fact that for any $s > 0$ and any $f,g\in\sw,$ 
both sides of \eqref{KPinfty} are finite, makes it quite difficult to find a counter-example to \eqref{KPinfty}. Such counter-example should violate the structure of the right-hand side of \eqref{KPinfty}, but not the fact that the left-hand side is finite. As a step towards \eqref{KPinfty} the purpose of this note is to prove the following results

\begin{thm}\label{thm:bGN} Let $0 \leq r < s < t$ and set
\begin{equation}
\alpha := \frac{t - s}{t - r} \quad \text{ and } \quad \beta:= \frac{s-r}{t -r}.
\end{equation}
Then, for every $f, g \in \sw$ we have
\begin{equation}\label{KPIBinfty}
\norm{D^s(fg)}{L^\infty}  \lesssim \norm{D^{r}f}{\Bzinf}^{\alpha}  \norm{D^{t}f}{\Bzinf}^{\beta} \norm{g}{L^\infty} +  \norm{f}{L^\infty} \norm{D^{r}g}{\Bzinf}^{\alpha} \norm{D^{t}g}{\Bzinf}^{\beta},
\end{equation}
where the implicit constant depends only on $r, s, t$, and dimension $n$. In particular,
\begin{equation}\label{KPIinfty}
\norm{D^s(fg)}{L^\infty}  \lesssim \norm{D^{r}f}{L^\infty}^{\alpha}  \norm{D^{t}f}{L^\infty}^{\beta} \norm{g}{L^\infty} + \norm{f}{L^\infty} \norm{D^{r}g}{L^\infty}^{\alpha} \norm{D^{t}g}{L^\infty}^{\beta}.
\end{equation}
\end{thm}

\begin{remark} Inequality \eqref{KPIinfty} can be regarded as a combination of Leibniz-rule and interpolation (or bilinear Gagliardo-Nirenberg) inequalities. Notice that \eqref{KPIinfty} is weaker than \eqref{KPinfty}. Indeed, given $0 \leq r < s < t$, by the linear Gagliardo-Nirenberg inequality (see, for instance, Theorem 2.44 in \cite{BCD}), we have
\begin{equation}\label{linearinterprst}
\norm{D^s f}{L^\infty}  \lesssim \norm{D^rf}{L^\infty}^{\frac{t-s}{t-r}}  \norm{D^tf}{L^\infty}^{\frac{s-r}{t-r}}, \quad \forall f \in \sw.
\end{equation}
Then, it follows that \eqref{KPinfty}, if true, would imply \eqref{KPIinfty}.
\end{remark}

\begin{thm}\label{thm:sp1p2} Suppose $s > 2n+1$. Let $1 < p_1, p_2 < \infty$ and $\eps > 0$ with $n/p:=(1/p_1 + 1/p_2) n < \eps < 1$. Then for every $f, g \in \sw$ we have
\begin{align*}
\norm{D^s(fg)}{L^\infty}& \lesssim  \norm{D^{s}f}{L^{p_1}}^{1-\frac{n}{p \eps}} \norm{D^{s+\eps}f}{L^{p_1}}^{\frac{n}{p \eps}}\norm{g}{L^{p_2}}  +  \norm{f}{L^{p_1}} \norm{D^{s}g}{L^{p_2}}^{1-\frac{n}{p \eps}} \norm{D^{s+\eps}g}{L^{p_2}}^{\frac{n}{p \eps}} \\
& + \norm{D^sf}{L^\infty} \norm{g}{L^\infty} +  \norm{f}{L^\infty} \norm{D^sg}{L^\infty},
\end{align*}
where the implicit constant depends only on $s$, $n$, $\eps$, $p_1$, and $p_2$.
\end{thm}

\begin{remark} In the case $s > 2n +1$, the proof of Theorem \ref{thm:sp1p2} will be based on a connection between Kato-Ponce inequalities and the bilinear Calder\'on-Zygmund theory, see Section \ref{secc:2n+1}. Notice that the inequality in Theorem \ref{thm:sp1p2} involves no derivatives lower than $D^s$. Also, $\eps >0$ can be arbitrarily small and $p_1, p_2 \in (1,\infty)$ arbitrarily large, as long as $(1/p_1 + 1/p_2) n < \eps$.
\end{remark}

\section{Preliminaries}

Let $\Phi: \rn \to \re$ be a smooth, non-negative, radial function supported in $\{\xi\in\re^n: |\xi| \leq 2\}$ with $\Phi \equiv 1$ in $\{\xi\in\re^n: |\xi| \leq 1\}.$ Define $\Psi: \rn \to \re$ supported in $1/2 \leq |\xi| \leq 2$ as $\Psi(\xi):=\Phi(\xi) - \Phi(2\xi)$ for  $\xi \in \rn$,  so that
\begin{equation}\label{AOTI}
\sum\limits_{j \in \ent} \Delta_j h = h \, \text{ in } \mathcal{S}'(\rn)  \quad \forall h \in \sw,
\end{equation}
where, as usual,  $\Delta_jh$ is defined for $h \in \sw$  as
$$
\widehat{\Delta_j h}(\xi) := \Psi(2^{-j}\xi) \widehat{h}(\xi) \quad \forall \xi \in \rn. 
$$

We recall that the  Besov $\Bzinf$-norm is given by
\begin{equation}\label{embeddings}
\norm{h}{\Bzinf}:=  \sup\limits_{j \in \ent} \norm{\Delta_j h}{L^\infty} \leq  \norm{\widehat{\Psi}}{L^1} \norm{h}{L^\infty}.
\end{equation}
For $f \in \sw$ and $\lambda > 0$ set $f_\lambda(x):=f(\lambda x)$ for every $x \in \rn$.  For $s \geq 0$ we have 
\begin{equation}\label{scaling}
  \norm{D^s(f_\lambda)}{\Bzinf} = \lambda^s \norm{D^sf}{\Bzinf}\quad \text{ for all } \lambda =2^{j_0},\,j_0\in\ent.
\end{equation}

We note tha $\tilde{\Phi}(\xi + \eta) \Phi(\xi) \Psi(\eta) = \Phi(\xi) \Psi(\eta)$ for every $\xi, \eta \in \rn,$ where $\tilde{\Phi}(\cdot):=\Phi(4^{-1}\cdot),$ and  
 write $\Phi_{(s)}(\cdot) := |\cdot|^s \tilde{\Phi}(\cdot) $. Reasoning as in \cite{GO}, the absolutely convergent Fourier series for $\Phi_{(s)}(t)\chi_{[-8,8]^n}(t),$ 
\begin{equation}\label{FourierPhi}
\Phi_{(s)}(t) = \sum_{m \in \zn} c_{s,m} e^{\frac{2\pi i}{16} m \cdot t}\chi_{[-8,8]^n}(t),
\end{equation}
has coefficients $c_{s,m}$ satisfying 
\begin{equation}\label{sizecm}
c_{s,m} = O(1+|m|^{-n-s}).
\end{equation}

\section{Proof of Theorem \ref{thm:bGN}}\label{secc:bGN}

\begin{proof} Fix $0 \leq r < s < t$.  By \eqref{AOTI}, we have
$$
D^s(fg)(x) = \int_{\re^{2n}} |\xi + \eta|^s \fh(\xi) \gh(\eta) e^{2 \pi i (\xi + \eta) \cdot x} d\xi d\eta =: \Pi(f,g)(x) + \tilde{\Pi}(f,g)(x),
$$
with
$$
\Pi(f,g)(x):=  \int_{\re^{2n}} \sum_{j \in \zz} \sum_{k \le j} |\xi + \eta|^s \Psi(2^{-j}\xi )\Psi(2^{-k} \eta)\fh(\xi) \gh(\eta) e^{2 \pi i (\xi + \eta) \cdot x} d\xi d\eta
$$
and
$$
 \tilde{\Pi}(f,g)(x):= \int_{\re^{2n}} \sum_{j \in \zz} \sum_{j < k} |\xi + \eta|^s \Psi(2^{-j}\xi )\Psi(2^{-k} \eta)\fh(\xi) \gh(\eta) e^{2 \pi i (\xi + \eta) \cdot x} d\xi d\eta.
$$
Now, we split $\Pi(f,g)$ (and then, similarly, $\tilde{\Pi}$) as follows
\begin{align*}
\Pi(f,g)(x) & =  \int_{\re^{2n}}  \sum_{j \in \zz} |\xi + \eta|^s \Psi(2^{-j}\xi )\Phi(2^{-j} \eta)\fh(\xi) \gh(\eta) e^{2 \pi i (\xi + \eta) \cdot x} d\xi d\eta\\
& =  \int_{\re^{2n}}  \sum_{j \leq 0}|\xi + \eta|^s \Psi(2^{-j}\xi )\Phi(2^{-j} \eta)\fh(\xi) \gh(\eta) e^{2 \pi i (\xi + \eta) \cdot x} d\xi d\eta\\
&+   \int_{\re^{2n}} \sum_{j >0} |\xi + \eta|^s \Psi(2^{-j}\xi )\Phi(2^{-j} \eta)\fh(\xi) \gh(\eta) e^{2 \pi i (\xi + \eta) \cdot x} d\xi d\eta\\
& =  \int_{\re^{2n}}  \sum_{j \leq 0} \frac{|\xi + \eta|^s}{|\xi|^r} \Psi(2^{-j}\xi )\Phi(2^{-j} \eta)\widehat{D^rf}(\xi) \gh(\eta) e^{2 \pi i (\xi + \eta) \cdot x} d\xi d\eta\\
&+  \int_{\re^{2n}}  \sum_{j > 0} \frac{|\xi + \eta|^s}{|\xi|^t} \Psi(2^{-j}\xi )\Phi(2^{-j} \eta)\widehat{D^tf}(\xi) \gh(\eta) e^{2 \pi i (\xi + \eta) \cdot x} d\xi d\eta\\
& =: \Pi_1(D^rf, g) + \Pi_2(D^tf,g).
\end{align*}
We now look at the bilinear kernel of $\Pi_1$ (the kernel for $\Pi_2$ will be dealt with in a similar way).
\begin{equation}\label{defPi1}
\Pi_1(f,g)(x) = \int_{\re^{2n}} K_1(x-y,x-z) f(y) g(z) dy dz,
\end{equation}
where, after putting $\Psi_{(-r)}(\cdot):= |\cdot|^{-r} \Psi(\cdot)$ and using that $\tilde{\Phi}(\xi+\eta) \Phi(\xi) \Psi(\eta)=\Phi(\xi) \Psi(\eta)$ for every $\xi, \eta \in \rn$,  
 $K_1$ is given by
\begin{align*}
K_1(y,z) & = \int_{\re^{2n}} \sum_{j \leq 0} \frac{|\xi + \eta|^s}{|\xi|^r} \Psi(2^{-j}\xi )\Phi(2^{-j} \eta) e^{2 \pi i (\xi  \cdot y + \eta \cdot z)} d\xi d\eta\\
& = \int_{\re^{2n}} \sum_{j \leq 0} \frac{2^{js}}{2^{jr}} \Phi_{(s)}(2^{-j}(\xi + \eta)) {\Psi_{(-r)}}(2^{-j}\xi )\Phi(2^{-j} \eta) e^{2 \pi i (\xi  \cdot y + \eta \cdot z)} d\xi d\eta.
\end{align*}
Hence, using the  Fourier expansion in \eqref{FourierPhi} and noting that the support of  $\psi_{(-r)}(\xi)\phi(\eta)$ is contained in $\{(\xi,\eta):|\xi+\eta|\le 4\},$ we get
\begin{align*} 
K_1(y,z) & =  \int_{\re^{2n}} \sum_{j \leq 0} \sum_{m \in \zn } c_{s,m} 2^{j(s-r)} e^{\frac{2\pi i}{16}  m \cdot 2^{-j}(\xi + \eta)} \Psi_{(-r)}(2^{-j}\xi )\Phi(2^{-j} \eta) e^{2 \pi i (\xi  \cdot y + \eta \cdot z)} d\xi d\eta\\
& =   \sum_{j \leq 0} \sum_{m \in \zn } c_{s,m} 2^{j(s-r)} \int_{\re^{2n}} e^{\frac{2\pi i}{16}  m \cdot 2^{-j}(\xi + \eta)} \Psi_{(-r)}(2^{-j}\xi )\Phi(2^{-j} \eta) e^{2 \pi i (\xi  \cdot y + \eta \cdot z)} d\xi d\eta\\
& =   \sum_{j \leq 0} \sum_{m \in \zn } c_{s,m} 2^{j(s-r)} 2^{2jn} \int_{\re^{2n}} e^{\frac{2\pi i}{16}  m \cdot (\xi + \eta)} \Psi_{(-r)}(\xi )\Phi(\eta) e^{2 \pi i 2^j (\xi  \cdot y + \eta \cdot z)} d\xi d\eta\\
& =   \sum_{j \leq 0} \sum_{m \in \zn } c_{s,m} 2^{j(s-r)} 2^{2jn} \widehat{\Psi_{(-r)}}(\fr{m}{16} + 2^j y) \widehat{\Phi}(\fr{m}{16} + 2^{j}z).
\end{align*}
Now,
\begin{align*}
\Pi_1(f,g)(x) & = \sum\limits_{l \in \ent} \int_{\re^{2n}} K_1(x-y,x-z) (\Delta_{l} f)(y) g(z) dy dz\\
& \leq  \sum_{j \leq 0} \sum_{m \in \zn }  \sum\limits_{l \in \ent }c_{s,m} 2^{j(s-r)} \\
&\times \int_{\re^{2n}} 2^{2jn}  \widehat{\Psi_{(-r)}}(\fr{m}{16} + 2^j (x-y)) \widehat{\Phi}(\fr{m}{16} + 2^{j}(x-z)) \Delta_{l} f(y) g(z) dy dz.
\end{align*}
For a fixed $j \in \ent$ we look at the integral in $y$
\begin{align*}
\int_{\rn}  \widehat{\Psi_{(-r)}}(\fr{m}{16} + 2^j (x-y)) \Delta_{l}f(y) dy = \int_{\rn} \frac{e^{2\pi i \xi \cdot (2^{-j} \fr{m}{16} + x)}}{2^{jn}}  \Psi(2^{-j} \xi) \Psi(2^{-l} \xi) \hat{{f}}(\xi) d\xi,
\end{align*}
which, due to the support conditions on $\Psi$, vanishes for every $l \in \ent \setminus \{j-1, j, j+1\}$. Consequently, 
\begin{align*}
& |\Pi_1(f,g)(x)| \leq  \sum_{j \leq 0} \sum_{m \in \zn }  \sum\limits_{l = j-1, j, j+1} |c_{s,m}| 2^{j(s-r)} \\
&\times  \int_{\re^{2n}} 2^{2jn}  |\widehat{\Psi_{(-r)}}(\fr{m}{16} + 2^j (x-y))| |\widehat{\Phi}(\fr{m}{16} + 2^{j}(x-z))| |\Delta_{l} f(y)| |g(z)| dy dz\\
& \leq  3 \left(\sum_{j \leq 0}  2^{j(s-r)} \right) \left(\sum_{m \in \zn } |c_{s,m}|\right) \norm{\widehat{\Psi_{(-r)}}}{L^1(\rn)} \norm{\widehat{\Phi}}{L^1(\rn)} \norm{f}{\dot{B}^{0,\infty}_\infty} \norm{g}{L^\infty}.
\end{align*}
Since $s - r > 0$ we have $\sum_{j \leq 0}  2^{j(s-r)} < \infty$ and, from \eqref{sizecm}, $\sum_{m \in \zn } |c_{s,m}| < \infty$. Hence,\begin{align*}
|\Pi_1(f,g)(x)| & \leq C \norm{\widehat{\Psi_{(-r)}}}{L^1(\rn)} \norm{\widehat{\Phi}}{L^1(\rn)} \norm{f}{\dot{B}^{0,\infty}_\infty} \norm{g}{L^\infty} \quad \forall x \in \rn,
\end{align*}
where $C > 0$ depends only on $r, s$, and  $n$.

Along the same lines, now for $s < t$  one gets the bound for $\Pi_2(f,g)$,
\begin{align*}
|\Pi_2(f,g)(x)| & \leq  c \left(\sum_{j > 0}  2^{j(s-t)} \right) \left(\sum_{m \in \zn } |c_{s,m}|\right) \norm{{\widehat{\Psi_{(-t)}}}}{L^1(\rn)} \norm{\widehat{\Phi}}{L^1(\rn)}  \norm{f}{\dot{B}^{0,\infty}_\infty}  \norm{g}{L^\infty},
\end{align*}
with $s - t < 0$. Then
\begin{equation}\label{eq:pi}
\norm{\Pi(f,g)}{L^\infty} \leq C ( \norm{D^rf}{\dot{B}^{0,\infty}_\infty} + \norm{D^tf}{\dot{B}^{0,\infty}_\infty}) \norm{g}{L^\infty}.
\end{equation}
Interchanging the roles of $f$ and $g$ to deal with $\tilde{\Pi}$ yields
\begin{equation}\label{eq:pitilde}
\norm{\tilde{\Pi}(f,g)}{L^\infty} \leq C ( \norm{D^rg}{\dot{B}^{0,\infty}_\infty} + \norm{D^tg}{\dot{B}^{0,\infty}_\infty}) \norm{f}{L^\infty}.
\end{equation}

Given a positive dyadic number $\mu$, plugging in $f_{\mu}$ and $g_{\mu}$ into \eqref{eq:pi} and \eqref{eq:pitilde}, using the scaling property \eqref{scaling} and the fact that $\Pi(f_\mu,g_{\mu})=\mu^s\,\Pi(f,g)_\mu$ and 
$\tilde{\Pi}(f_\mu,g_{\mu})=\mu^s\,\tilde{\Pi}(f,g)_\mu$, we get
\begin{align*}
\norm{\Pi(f,g)}{L^\infty}  &\lesssim (\lambda^{r-s} \norm{D^rf}{\Bzinf} + \lambda^{t-s} \norm{D^tf}{\Bzinf}) \norm{g}{L^\infty},\\
\norm{\tilde{\Pi}(f,g)}{L^\infty} &\lesssim ( \lambda^{r-s} \norm{D^rg}{\Bzinf} + \lambda^{t-s} \norm{D^tg}{\Bzinf})\norm{f}{L^\infty},
\end{align*}
for every positive number $\lambda.$
Minimizing in $\lambda$  each of the above inequalities leads to 
\begin{align*}
\norm{\Pi(f,g)}{L^\infty}  &\lesssim \norm{D^{r}f}{\Bzinf}^{\alpha}  \norm{D^{t}f}{\Bzinf}^{\beta} \norm{g}{L^\infty}, \\
\norm{\tilde{\Pi}(f,g)}{L^\infty} &\lesssim   \norm{f}{L^\infty} \norm{D^{r}g}{\Bzinf}^{\alpha} \norm{D^{t}g}{\Bzinf}^{\beta},
\end{align*}
from which  \eqref{KPIBinfty} follows.
\end{proof}

\section{The case $s > 2n+1$}\label{secc:2n+1}

A smooth function $\sigma : \re^{2n} \setminus \{(0,0)\} \rightarrow \com$ is said to belong to the class of bilinear Coifman-Meyer multipliers if for all multi-indices $\alpha, \beta \in \na_0^n$ with $|\alpha|+|\beta| \leq 2n+1$ there exist constants $c_{\alpha, \beta} > 0$ such that
\begin{equation}\label{defCM}
|\partial_\xi^\alpha \partial_\eta^\beta \sigma(\xi, \eta)|  \leq c_{\alpha, \beta} (|\xi|+|\eta|)^{- |\alpha | - |\beta|}, \quad \forall (\xi, \eta) \in \re^{2n} \setminus \{(0,0)\}.
\end{equation}
In \cite{GO}, the bilinear mapping $(f,g) \mapsto D^s(fg)$ was decomposed into the sum of three bilinear multipliers as follows
\begin{equation}\label{decomp}
D^s(fg) = T_{1,s}(D^s f,g) + T_{2,s}(f,D^s g) + T_{3,s}(f,D^sg),
\end{equation}
where, keeping with the notation in Section \ref{secc:bGN}, for $(\xi, \eta) \in \re^{2n} \setminus \{(0,0)\}$ the bilinear multipliers for $T_{1,s}$ and $T_{2,s}$ are given by
\begin{equation}\label{symbolsT}
\sigma_{1,s}(\xi, \eta):= \sum\limits_{j \in \ent} \Psi(2^{-j}\xi) \Phi(2^{-j+3} \eta) \frac{|\xi + \eta|^s}{|\xi|^s} \quad \text{and} \quad \sigma_{2,s}(\xi, \eta):=\sigma_{1,s}(\eta,\xi),
\end{equation}
respectively, which belong to the Coifman-Meyer class for every $s > 0$. On the other hand, the multiplier for $T_{3,s}$, denoted by $\sigma_{3,s}$, can be expressed as
\begin{equation}\label{sigma3s}
\sigma_{3,s}(\xi, \eta) := \sum_{k \in \zz} \sum_{m \in \zz^n} c_{s,m} e^{\frac{2 \pi i }{16} 2^{-k} (\xi + \eta) \cdot m} \Psi(2^{-k}\xi) \Psi_{(-s)}(2^{-k}\eta).
\end{equation}
For fixed $\xi, \eta \in \re^{2n} \setminus \{(0,0)\}$ the condition on the support of $\Psi$ implies that the sum in $k$ has only finitely many terms; namely, those with $2^k \sim |\xi| \sim |\eta|$. When derivatives in $\xi$ and $\eta$ of the product $e^{\frac{2 \pi i }{16} 2^{-k} (\xi + \eta) \cdot m} \Psi(2^{-k}\xi) \Psi_{(-s)}(2^{-k}\eta)$ are taken, after each derivative a factor $2^{-k} \,  (\sim |\xi|^{-1} \sim |\eta|^{-1} \sim (|\xi| + |\eta|)^{-1})$ appears, producing the right-hand side of \eqref{defCM}. However, when the derivatives fall on the factor $e^{\frac{2 \pi i }{16} 2^{-k} (\xi + \eta) \cdot m}$ also components of $m \in \zz^n$ appear. Since the definition of a Coifman-Meyer multiplier requires at most $2n+1$ derivatives, the worst case scenario for the sum over $m \in \zz^n$ (i.e., the case in which all $2n+1$ derivatives fall on $e^{\frac{2 \pi i }{16} 2^{-k} (\xi + \eta) \cdot m}$) leads to the sum
$$
\sum_{m \in \zz^n} |c_{s,m}| |m|^{2n+1}.
$$
By \eqref{sizecm}, the sum above will be finite provided that $s > 2n + 1$. That is, whenever $s > 2n+1$ all three bilinear operators in \eqref{decomp}, and therefore the mapping $(f,g) \mapsto D^s(fg)$, can be realized as bilinear Coifman-Meyer multipliers. Since the class of Coifman-Meyer multipliers is included in the family of bilinear Calder\'on-Zygmund operators (see, \cite[Section 6]{GT}) all the mapping properties of the type
\begin{equation}\label{mapT}
\norm{T(f,g)}{Z} \lesssim \norm{f}{X} \norm{g}{Y},
\end{equation}
that apply to bilinear C-Z operators $T$ on function spaces $X, Y,$ and $Z$ will also apply to $(f,g) \mapsto D^s(fg)$. For example, for a bilinear C-Z operator $T$, given $ 1 < p_1, p_2 < \infty $ and $1/p:=1/{p_1} + 1/{p_2}$, it holds that
\begin{equation}\label{bilinearLp}
\norm{T(f,g)}{L^{p}}  \lesssim \norm{f}{L^{p_1}} \norm{g}{L^{p_2}}
\end{equation}
and (see \cite[Proposition 1]{GT}) that, 
\begin{equation}\label{infinfbmo}
\norm{T(f,g)}{BMO} \lesssim \norm{f}{L^\infty} \norm{g}{L^\infty},
\end{equation}
as well as other end-point estimates such as
\begin{equation}
\norm{T(f,g)}{L^{1,\infty}}  \lesssim \norm{D^{s}f}{L^\infty} \norm{g}{L^1} + \norm{f}{L^1} \norm{D^{s}g}{L^\infty}.
\end{equation}

As a consequence of the results above, we have

\begin{thm}\label{thm:bmo} If $s > 2n+1$, then for every $f, g \in \sw$ we have the endpoint inequalities
\begin{equation}\label{KPbmo}
\norm{D^s(fg)}{BMO}  \lesssim \norm{D^{s}f}{L^\infty} \norm{g}{L^\infty} + \norm{f}{L^\infty} \norm{D^{s}g}{L^\infty}.
\end{equation}
and
\begin{equation}
\norm{D^s(fg)}{L^{1,\infty}}  \lesssim \norm{D^{s}f}{L^\infty} \norm{g}{L^1} + \norm{f}{L^1} \norm{D^{s}g}{L^\infty}.
\end{equation}
\end{thm}

\begin{remark}\label{rmkonCZ}
We note that the conditions \eqref{defCM} being satisfied with up to $n+1$ derivatives (instead of $2n+1$) are sufficient for the corresponding multiplier operator to be bounded from $L^{p_1}\times L^{p_2}$ into $L^p$ for $1<p_1,p_2,p<\infty$ and $\frac{1}{p}=\frac{1}{p_1}+\frac{1}{p_2},$ as shown in Tomita~\cite{T}. The endpoint boundedness $L^\infty\times L^\infty$ into $BMO$ for Coifman-Meyer multipliers, with only up to $n+1$ derivatives in \eqref{defCM}, is unknown to us. To pass through the bilinear C-Z theory, as done above, it suffices that the conditions \eqref{defCM} be satisfied with up to $2n+1$ derivatives.
\end{remark}

\emph{Proof of Theorem \ref{thm:sp1p2}.} By hypothesis, $1/p:=1/p_1 + 1/p_2$, so that $n/p < \eps < 1$. It was proved in \cite[pp.193--198]{KoTa2}  that a function $F$ with $ \norm{D^\eps F}{L^p} +  \norm{F}{BMO} + \norm{F}{L^p} < \infty$ can be written as $F = F_0 + G + F_1$ where
\begin{equation}\label{boundKoTa2}
\norm{F_0}{L^\infty} \lesssim \norm{D^\eps F}{L^p},\quad \norm{G}{L^\infty} \lesssim \norm{F}{BMO}, \quad \text{and} \quad \norm{F_1}{L^\infty} \lesssim \norm{F}{L^p}.
\end{equation}
Now, with $T_{1,s}$ as in the decomposition \eqref{decomp}, let us first choose $F:=T_{1,s}(D^s f, g)$, so that from \eqref{boundKoTa2} we get
$$
\norm{T_{1,s}(D^s f,g)}{L^\infty} \lesssim \norm{T_{1,s}(D^sf,g)}{L^p} +  \norm{T_{1,s}(D^sf,g)}{BMO} +  \norm{D^{\eps}(T_{1,s}(D^sf,g))}{L^p}.
$$
The fact that $T_{1,s}$ is a bilinear C-Z operator and \eqref{bilinearLp} yield
$$
 \norm{T_{1,s}(D^sf,g)}{L^p}  \lesssim \norm{D^{s}f}{L^{p_1}} \norm{g}{L^{p_2}}.
$$
Also, from \eqref{infinfbmo}, it follows that
$$
\norm{T_{1,s}(D^sf,g)}{BMO} \lesssim  \norm{D^sf}{L^\infty} \norm{g}{L^\infty}.
$$
On the other hand, notice that
$$
D^\eps(T_{1,s}(D^s f,g)) =: T_{1, s+\eps} (D^{s+\eps}f,g),
$$
where the bilinear symbol for the operator $T_{1, s+\eps}$ equals $\sigma_{1, s + \eps}(\xi, \eta)$ (using the notation in \eqref{symbolsT}), also a Coifman-Meyer multiplier. Hence, \eqref{bilinearLp} gives
$$
 \norm{D^\eps(T_{1,s}(D^s f,g))}{L^p}  \lesssim \norm{D^{s+\eps}f}{L^{p_1}} \norm{g}{L^{p_2}}.
$$
Putting all together, for $T_{1,s}(f,g)$ we have
\begin{equation}\label{boundT1s0}
\norm{T_{1,s}(D^s f,g)}{L^\infty} \lesssim  (\norm{D^{s}f}{L^{p_1}} + \norm{D^{s+\eps}f}{L^{p_1}}) \norm{g}{L^{p_2}} + \norm{D^sf}{L^\infty} \norm{g}{L^\infty}.
\end{equation}
Given a positive dyadic number $\mu$, by replacing $f$ and $g$ in \eqref{boundT1s0} with $f_\mu$ and $g_\mu$ and using the facts that
$$
\norm{D^s(f_\mu)}{L^q} =\mu^{s-\frac{n}{q}} \norm{D^s f}{L^q}, \quad \forall q \in [1,\infty],
$$ 
that $1/p=1/p_1+1/p_2$, and that $T_{1,s}(D^sf_\mu,g_{\mu})=\mu^s\,T_{1,s}(D^sf,g)_\mu,$  we obtain
$$
\norm{T_{1,s}(D^s f,g)}{L^\infty} \lesssim  (\lambda^{-\frac{n}{p}} \norm{D^{s}f}{L^{p_1}} + \lambda^{\eps - \frac{n}{p}}\norm{D^{s+\eps}f}{L^{p_1}}) \norm{g}{L^{p_2}} + \norm{D^sf}{L^\infty} \norm{g}{L^\infty},
$$
for every positive number $\lambda.$
Minimization over $\lambda$ then implies 
\begin{equation}\label{boundT1s}
\norm{T_{1,s}(D^s f,g)}{L^\infty} \lesssim  \norm{D^{s}f}{L^{p_1}}^{1-\frac{n}{p \eps}} \norm{D^{s+\eps}f}{L^{p_1}}^{\frac{n}{p \eps}}\norm{g}{L^{p_2}} + \norm{D^sf}{L^\infty} \norm{g}{L^\infty}.
\end{equation}
And, by an analogous argument based on $T_{2,s}$,
\begin{equation}\label{boundT2s}
\norm{T_{2,s}(f,D^s g)}{L^\infty} \lesssim  \norm{f}{L^{p_1}} \norm{D^{s}g}{L^{p_2}}^{1-\frac{n}{p \eps}} \norm{D^{s+\eps}g}{L^{p_2}}^{\frac{n}{p \eps}}  +  \norm{f}{L^\infty} \norm{D^sg}{L^\infty}.
\end{equation}
It only remains to consider $T_{3,s}$. Since $s > 2n +1$, again from \eqref{bilinearLp} and \eqref{infinfbmo}, we have
$$
\norm{T_{3,s}(D^s f,g)}{L^p} +  \norm{T_{3,s}(D^s f,g)}{BMO} \lesssim  \norm{D^{s}f}{L^{p_1}} \norm{g}{L^{p_2}} +  \norm{D^sf}{L^\infty} \norm{g}{L^\infty}.
$$
Now,  
$$
D^\eps(T_{3,s}(D^s f,g)) =: T_{3, s+\eps} (D^{s+\eps}f,g)
$$
where the bilinear symbol for $T_{3, s+\eps}$ is similar to $\sigma_{3,s}$ in \eqref{sigma3s} but with $c_{s,m}$ replaced by $c_{s+\eps, m}$, the Fourier coefficients for $\Phi_{(s+\eps)}$ which will satisfy $c_{s+\eps,m} = O(1+|m|^{-n-s-\eps})$. Consequently, 
$$
 \norm{D^\eps(T_{3,s}(D^s f,g))}{L^p}  \lesssim \norm{D^{s+\eps}f}{L^{p_1}} \norm{g}{L^{p_2}}
$$
and, proceeding as before, after scaling we get
\begin{equation}\label{boundT3s}
\norm{T_{3,s}(D^s f,g)}{L^\infty} \lesssim  \norm{D^{s}f}{L^{p_1}}^{1-\frac{n}{p \eps}} \norm{D^{s+\eps}f}{L^{p_1}}^{\frac{n}{p \eps}}\norm{g}{L^{p_2}} + \norm{D^sf}{L^\infty} \norm{g}{L^\infty}.
\end{equation}
Finally, Theorem \ref{thm:sp1p2} follows from \eqref{decomp}, \eqref{boundT1s}, \eqref{boundT2s}, and \eqref{boundT3s}.\qed

\section*{Acknowledgements} The authors are grateful to Fr\'ed\'eric Bernicot and Gustavo Ponce for useful conversations.

\bibliographystyle{amsplain}

\begin{thebibliography}{99}

\bibitem{BBR} N. Badr, F. Bernicot, and E. Russ, \emph{Algebra properties for Sobolev spaces - Applications to semilinear PDE's on manifolds}, J. Anal. Math. {\bf{118}} (2012), 509-544. 

\bibitem{BCD} H. Bahouri, J.-Y. Chemin, and R. Danchin. Fourier Analysis and Nonlinear Partial Differential Equations. Grundlehren der mathematischen Wissenschaften, Volume 343. Springer-Verlag, 2011. 

\bibitem{BMMN}  F. Bernicot, D. Maldonado, K. Moen, and V. Naibo, \emph{Bilinear Sobolev-Poincar\'e inequalities and Leibniz-type rules}, to appear in J. Geom. Anal., DOI: 10.1007/s12220-012-9367-4.

\bibitem{CWein} M. Christ and M. Weinstein,
\newblock {\it Dispersion of small amplitude solutions of the generalized Korteweg-de Vries equation},
\newblock {J. Funct. Anal.} {\bf 100}, (1991), 87--109.


\bibitem{CRT} T. Coulhon, E. Russ, and V. Tardivel-Nachef, \emph{Sobolev algebras on Lie groups and Riemannian manifolds}, Amer. J. Math., {\bf{123}}, (2001), 283--342.



\bibitem{GO} L. Grafakos and S. Oh, \emph{The Kato-Ponce inequality,} to appear in Comm. PDE.,  DOI:10.1080/03605302.2013.822885.

\bibitem{GT} L.\ Grafakos and R.\ H.\ Torres, \emph{Multilinear Calder\'{o}n-Zygmund theory,} Adv.\ in \ Math. {\bf 165} (2002), 124-164.

\bibitem{GK} A. Gulisashvili and M. Kon, \emph{Exact smoothing properties of Schr\"odinger semigroups}, Amer. J. Math. {\bf{118}}, (1996), 1215--1248.

\bibitem{KaPo} T. Kato and G. Ponce,
\newblock {\it Commutator estimates and the Euler and Navier-Stokes equations.}
\newblock {Comm. Pure Appl. Math.} {\bf 41}, (1988),  891--907.


\bibitem{KPVe}
C. Kenig, G. Ponce, and L. Vega,
\newblock {\it Well-posedness and scattering results for the generalized Korteweg-de Vries equation via the contraction principle},
\newblock {Comm. Pure Appl. Math.} {\bf 46} (1993), 527--620.

\bibitem{KoTa1} H. Kozono and Y. Taniuchi, \emph{Bilinear estimates in BMO and the Navier-Stokes equations},  Math. Z. {\bf{235}}, (2000), 173--194.

\bibitem{KoTa2} H. Kozono and Y. Taniuchi, \emph{Limiting case of the Sobolev inequality in BMO, with application to the Euler equations}, Commun. Math. Phys. {\bf{214}}, (2000),191--200.


\bibitem{T}
N. Tomita,
\newblock {\it A H\"ormander type multiplier theorem for multilinear operators},
\newblock {J. Funct. Anal.} {\bf 259} (2010), no. 8, 2028--2044.



\end{thebibliography}

\end{document}